\documentclass[a4paper,11pt,reqno,twoside]{amsart}
\usepackage{amsmath}
\usepackage{amsfonts}
\usepackage{amssymb}
\usepackage{amsthm}

\theoremstyle{plain}
\newtheorem{thm}{Theorem}[section]

\newtheorem{lemma}[thm]{Lemma}
\newtheorem{prop}[thm]{Proposition}
\newtheorem{defn}[thm]{Definition}
\theoremstyle{definition}
\newtheorem{remark}[thm]{Remark}
{\theoremstyle{definition} \newtheorem*{example}{Example}
                        \newtheorem*{claim}{Claim}}

\numberwithin{equation}{section}

\newcommand{\C}{\mathbb{C}}                         
\newcommand{\Z}{\mathbb{Z}}                         
\newcommand{\Q}{\mathbb{Q}}
\newcommand{\h}{\mathfrak{h}}                       
\newcommand{\F}{\mathcal{F}}
\renewcommand{\det}{\mathrm{det}}                   

\newcounter{mylist}

\newenvironment{ellbin}[2]
{\begin{bmatrix} #1\\ #2\end{bmatrix}}

\begin{document}
\date{\today}
\title[Elliptic $U(2)$ quantum group and elliptic hypergeometric series]
{Elliptic $U(2)$ quantum group and \\ elliptic hypergeometric series}
\author{Erik Koelink}
\author{Yvette van Norden}
\thanks{The second author is supported by Netherlands Organisation for Scientific
Research (NWO) under project number 613.006.572.}
\address{Technische Universiteit Delft, Faculteit Elektrotechniek, Wiskunde en Informatica,
 Toegepaste Wis\-kun\-dige Analyse,
Postbus 5031, 2600 GA Delft, the Netherlands}
\email{h.t.koelink@math.tudelft.nl, y.vannorden@math.tudelft.nl}

\author{Hjalmar Rosengren}
\address{Department of Mathematics, Chalmers University of
Technology and G{\"o}teborg University, SE-412 96 G{\"o}teborg,
Sweden} \email{hjalmar@math.chalmers.se}
\begin{abstract}
We investigate an elliptic quantum group introduced by Felder and
Varchenko, which is constructed from the $R$-matrix of the
Andrews--Baxter--Forrester model, containing both spectral and
dynamical parameter. We explicitly compute the matrix elements of
certain corepresentations and obtain orthogonality relations for
these elements. Using dynamical representations these
ortho\-go\-na\-lity relations give discrete bi-orthogonality
relations for terminating very-well-poised balanced elliptic
hypergeometric series, previously obtained by Frenkel and Turaev
and by Spiridonov and Zhedanov in different contexts.
\end{abstract}
\maketitle
\section{Introduction}

Elliptic functions appear in various solvable models in
statistical mechanics and other areas of physics. A famous example
is Baxter's $8$-vertex model \cite{b}, whose $R$-matrix,
containing the Boltzmann weights, is an elliptic solution of the
Yang--Baxter equation. A related face model was introduced by
Andrews, Baxter and Forrester \cite{abf}. In this case the
$R$-matrix satisfies a modified, ``dynamical'', version of the
Yang--Baxter equation,
 generalizing Wigner's hexagon identity for the classical $6j$-symbols of
quantum mechanics.

In the early 1980's, the algebraic study of the Yang--Baxter
equation lead to the introduction of quantum groups. The most well
understood quantum groups are those constructed from the simplest,
constant, solutions.  Quantum groups connected to more complicated
solutions, and in particular to elliptic solutions, have been more
difficult to construct and study. One reason for this is that
elliptic quantum groups
 are not Hopf algebras. Various approaches
have been tried for finding a substitute; cf.
\cite{ef,fe,fo,JimboKOS,s}. In the dynamical case, a decisive step
was taken by Felder and Varchenko \cite{FelderVarchenko}, who
introduced  the algebra that we will study here.  This example
motivated Etingof and Varchenko \cite{EtingofVarchenko} to
introduce $\h$-Hopf algebroids, a generalization of Hopf algebras
adapted to studying  dynamical $R$-matrices; cf.\
\cite{KoelinkRosengren,r} for  further additions to this
framework.

An important mathematical application of quantum groups is their
relation to basic hypergeometric series (or $q$-series), a class
of special functions  going back to work of Cauchy and Heine in
the 1840's. The input from quantum group theory has been important
for the rapid development of this field during the last 20 years.
To our knowledge, nobody has so far associated special functions
to elliptic quantum groups in an analogous way. There is, however,
a natural candidate for the special functions that should appear,
namely, the elliptic or modular hypergeometric series of Frenkel
and Turaev \cite{FrenkelTuraev}. This type of sums may be used to
express the elliptic $6j$-symbols of Date et al.\ \cite{d}, which
 are solutions to the Yang--Baxter equation that greatly generalize the
Andrews--Baxter--Forrester solution. For more information on
elliptic hypergeometric series we refer to
 \cite{kn,r2,Rosengren,sp1,sp2,sp3,sp4,SpiriZhe,Warnaar}.

Our main aim is to give an explicit link from elliptic quantum
groups to elliptic hypergeometric series. Namely, we show that
${}_{10}\omega_9$ sums, or elliptic $6j$-symbols, appear as matrix
elements for an elliptic quantum group which we denote by
$\F_R(U(2))$, which is the algebra of Felder and Varchenko with
some extra structure. To achieve this, we first construct
finite-dimensional corepresentations  of  $\F_R(U(2))$, analogous
to the standard representations of $SU(2)$ on spaces of
homogeneous polynomials in two variables. A main result, Theorem
\ref{thm:corepmatrix}, is an explicit expression for the matrix
elements of these corepresentations.  We can then calculate the
action of the matrix elements in  representations found by Felder
and Varchenko, and show that it is given in terms of elliptic
hypergeometric series.

The matrix elements satisfy orthogonality relations in the
non-commutative algebra $\F_R(U(2))$. Evaluating these in a
representation leads to bi-orthogonality relations for
${}_{10}\omega_9$ series. These relations were found already by
Frenkel and Turaev \cite{FrenkelTuraev}; cf.\ also
\cite{SpiriZhe}.

Our new derivation of the bi-orthogonality relations shows that
they can be viewed  as analogues of the orthogonality relations
for Krawtchouk polynomials, see \cite{VilenkinKlimyk} for the Lie
group $SU(2)$. For the quantum $SU(2)$ group the same approach
leads to quantum $q$-Krawtchouk polynomials, see
\cite{Koornwinder}. For the dynamical quantum $SU(2)$ group, i.e.\
corresponding to a trigonometric dynamical $R$-matrix, we get the
orthogonality relations for  $q$-Racah polynomials, see \cite[\S
4] {KoelinkRosengren}. So the above cases can be considered as
limiting cases of the bi-orthogonality relations for  elliptic
$6j$-symbols.

The paper is organized as follows. In section $2$ we recall the
definition of an $\h$-Hopf algebroid and the generalized
FRST-construction from \cite{EtingofVarchenko}. Then we describe
the elliptic  quantum group $\F_R(U(2))$, which is obtained from
the $R$-matrix of the Andrews--Baxter--Forrester model.  In
section \nolinebreak[4] $3$ we define finite-dimensional
corepresentations of $\F_R(U(2))$ and compute their matrix
elements explicitly. In section $4$ we consider representations of
$\F_R(U(2))$ , from which we obtain commutative versions of  the
orthogonality relations for  matrix elements of the
corepresentations. It turns out that these are in fact
bi-orthogonality relations for terminating very-well-poised
balanced elliptic hypergeometric $_{10}\omega_{9}$-series (or
elliptic $6j$-symbols).

{\bf Notation:} We denote by $\theta(z)$ the normalized Jacobi
theta function
$$\theta(z)=\prod_{j=0}^\infty\left(1-zp^j\right)\left(1-p^{j+1}/z
\right),\qquad |p|<1,$$ where $p$ is a fixed parameter that is
suppressed from the notation. It satisfies
$$\theta(pz)=\theta(z^{-1})=-z^{-1}\theta(z),$$
 and the addition formula
    \begin{equation}\label{eqn:theta}
\theta(xy, x/y, zw, z/w)= \theta(xw,x/w,zy,z/y) +
(z/y)\theta(xz,x/z,yw,y/w),
    \end{equation}
 where we use the notation
$$\theta(a_1, \ldots ,a_n)= \theta(a_1)\cdots\theta(a_n).$$ We
define elliptic Pochhammer symbols by
$$(a)_n=\prod_{i=0}^{n-1}\theta(aq^{2i}), $$
 with $q$  another fixed parameter. We will frequently write
  $$(a_1,a_2,\dots,a_k)_n=
(a_1)_n\cdots(a_k)_n.$$ Elliptic binomial coefficients are defined
by
 $$\ellbin{k}{l}=\prod_{i=1}^l
\frac{\theta(q^{2(k-l+i)})}{\theta(q^{2i})}.$$
 Finally, the balanced very-well-poised elliptic hypergeometric series is
defined by \cite{FrenkelTuraev}
    \begin{equation*}
_{r+1}\omega_r(a_1; a_4, a_5,\ldots,a_{r+1}) =\sum_{k=0}^\infty
\frac{\theta(a_1 q^{4k})}{\theta(a_1)} \frac{(a_1, a_4, \ldots,
a_{r+1})_k q^{2k}}{(q^2,a_1q^2/a_4,\ldots,a_1q^2/a_{r+1})_k},
    \end{equation*}
where $(a_4\cdots a_{r+1})^2=a_1^{r-3}q^{2(r-5)}$. In this paper
all series terminate, i.e.\ one of the $a_i$ is of the form
$q^{-2n}$ with $n$ a nonnegative integer, so there are no
convergence problems. Let us emphasize that in this paper all
elliptic factorials and elliptic hypergeometric series are in base
$q^2,p$.

\section{Elliptic $U(2)$ quantum group}

In this section we recall the definition of $\h$-Hopf algebroids
(also known as dynamical quantum groups) and the
FRST-construction. We start with the definition of the quantum
dynamical Yang-Baxter equation with spectral parameter and give
in (\ref{eqn:ellR}) the $R$-matrix to which we apply the
FRST-construction. The generators and relations for the resulting
$\h$-Hopf algebroid have been studied by Felder and Varchenko
\cite{FelderVarchenko}. The paper \cite{FelderVarchenko} contains
hardly any proofs, so we provide a proof of one of their results
in Lemma \ref{lem:det}.

Let $\h$ be a finite dimensional complex vector space, viewed as a
commutative Lie algebra and $V=\bigoplus_{\alpha\in\h^*} V_\alpha$
a diagonalizable $\h$-module. The quantum dynamical Yang-Baxter
equation with spectral parameter (QDYBE) is given by
    \begin{equation}
    \begin{split}\label{eqn:DYBEs}
R^{12}(\lambda-h^{(3)}, z_{12})R^{13}(\lambda, z_{13})
R^{23}(\lambda-h^{(1)}, z_{23}) =R^{23}(\lambda,
z_{23})R^{13}(\lambda-h^{(2)}, z_{13}) R^{12}(\lambda, z_{12}).
    \end{split}
    \end{equation}
Here $R: \h^*\times \C\to \mathrm{End}(V\otimes V)$ is a
meromorphic function, $h$ indicates the action of $\h$, the upper
indices are leg-numbering notation for the tensor product and
$z_{ij}=z_i/z_j$. For instance, $R^{12}(\lambda-h^{(3)}, z)$
denotes the operator $R^{12}(\lambda-h^{(3)}, z) (u\otimes v
\otimes w )= R(\lambda-\mu, z)(u\otimes v) \otimes w$ for $w\in
V_\mu$. An $R$-matrix is by definition a solution of the QDYBE
(\ref{eqn:DYBEs}) which is $\h$-invariant.

In the example we study, $\h$ is one-dimensional. We identify
$\h=\h^*=\C$ and take $V$ to be the two-dimensional $\h$-module
$V=\C e_1\oplus \C e_{-1}$. In the basis $e_1\otimes e_1$,
$e_1\otimes e_{-1}$, $e_{-1}\otimes e_1$, $e_{-1}\otimes e_{-1}$
the $R$-matrix is given by
    \begin{equation}\label{eqn:ellR}
R(\lambda,z)=R(\lambda,z,p,q)=
\begin{pmatrix}
1 & 0 & 0 & 0 \\
 0 & a(\lambda,z) & b(\lambda,z) & 0 \\
0 & c(\lambda,z) & d(\lambda,z) & 0 \\
 0 & 0 & 0 & 1
\end{pmatrix},
    \end{equation}
where
    \begin{equation}
    \begin{split}\label{eqn:ellRelem}
a(\lambda,z)=\frac{\theta(z)\theta(q^{2(\lambda+2)})}
{\theta(q^2z)\theta(q^{2(\lambda+1)})},&\qquad
 b(\lambda,z)=\frac{\theta(q^2)\theta(q^{-2(\lambda+1)}z)}
 {\theta(q^2z)\theta(q^{-2(\lambda+1)})},\\
 c(\lambda,z)=\frac{\theta(q^2)\theta(q^{2(\lambda+1)}z)}
 {\theta(q^2z)\theta(q^{2(\lambda+1)})},&\qquad
 d(\lambda,z)=\frac{\theta(z)\theta(q^{-2\lambda})}{\theta(q^2z)\theta(q^{-2(\lambda+1)})}.
    \end{split}
    \end{equation}
The $R$-matrix defined by (\ref{eqn:ellR}) satisfies the QDYBE
(\ref{eqn:DYBEs}), see \cite{ba,abf,FelderVarchenko,JimboKOS}.

\subsection{$\h$-Hopf algebroids}

In this subsection we recall the notion of $\h$-bialgebroids and
$\h$-Hopf algebroids (or dynamical quantum groups) originally
introduced by Etingof and Varchenko \cite{EtingofVarchenko}, see
also \cite{EtingofSchiffmann}. For the definition of the antipode
in an $\h$-Hopf algebroid we follow \cite{KoelinkRosengren}. We
discuss the FRST-construction that associates an $\h$-bialgebroid
to an $\h$-invariant matrix, see \cite{EtingofVarchenko},
\cite{EtingofSchiffmann}.

\vspace{\baselineskip}

Let $\h$ be a finite dimensional complex vector space. Denote the
field of meromorphic functions on the dual of $\h$ by $M_{\h^*}$.

\begin{defn}
An $\h$-algebra is a complex associative algebra $A$ with $1$,
which is bigraded over $\h^*$, $A=\oplus_{\alpha,\beta\in\h^*}
A_{\alpha\beta}$, and equipped with two algebra embeddings
$\mu_l$, $\mu_r: M_{\h^*} \to A_{00}$ (the left and right moment
map) such that
    \begin{equation*}
\mu_l(f)a=a\mu_l(T_\alpha f), \qquad \mu_r(f)a= a\mu_r(T_\beta f),
\; \mbox{for all } a\in A_{\alpha\beta}, \; f\in M_{\h^*},
    \end{equation*}
where $T_\alpha$ denotes the automorphism $(T_\alpha f)(\lambda)=
f(\lambda+\alpha)$.

A morphism of $\h$-algebras is an algebra homomorphism preserving
the moment maps.
\end{defn}

Let $A$ and $B$ be two $\h$-algebras. The matrix tensor product
$A\tilde \otimes B$ is the $\h^*$-bigraded vector space with $(A
\tilde \otimes B)_{\alpha\beta}= \bigoplus_{\gamma\in\h^*}
(A_{\alpha\gamma}\otimes_{M_{\h^*}}B_{\gamma\beta})$, where
$\otimes_{M_{\h^*}}$denotes the usual tensor product modulo the
relations
    \begin{equation}\label{eqn:tensor}
\mu_r^A(f)a\otimes b= a\otimes \mu_l^B(f)b, \;\mbox{for all } a\in
A, b\in B, f \in M_{\h^*}.
    \end{equation}
The multiplication $(a\otimes b)(c\otimes d)=ac\otimes bd$ for
$a$, $c\in A$ and $b$, $d\in B$ and the moment maps
$\mu_l(f)=\mu_l^A(f)\otimes 1$ and $\mu_r(f)=1\otimes \mu_r^B(f)$
make $A\tilde \otimes B$ into an $\h$-algebra.

\begin{example}
Let $D_{\h}$ be the algebra of difference operators in $M_{\h^*}$,
consisting of the operators $\sum_i f_i T_{\beta_i}$, with $f_i\in
M_{\h^*}$ and $\beta_i\in\h^*$. This is an $\h$-algebra with the
bigrading defined by $fT_{-\beta} \in (D_{\h})_{\beta\beta}$ and
both moment maps equal to the natural embedding.

For any $\h$-algebra $A$, there are canonical isomorphisms $A
\cong A\tilde\otimes D_{\h} \cong D_{\h}\tilde\otimes A$, defined
by
    \begin{equation}\label{eqn:identDh}
x\cong x \otimes T_{-\beta} \cong T_{-\alpha}\otimes x, \;
\mbox{for all }x\in A_{\alpha\beta}.
    \end{equation}
The algebra $D_{\h}$ plays the role of the unit object in the
category of $\h$-algebras.
\end{example}

\begin{defn}
An $\h$-bialgebroid is an $\h$-algebra $A$ equipped with two
$\h$-algebra homomorphisms $\Delta:A\to A\tilde\otimes A$ (the
comultiplication) and $\varepsilon:A \to D_{\h}$ (the counit) such
that $(\Delta\otimes id)\circ\Delta= (id\otimes
\Delta)\circ\Delta$ and $(\varepsilon\otimes id)\circ
\Delta=id=(id\otimes \varepsilon)\circ\Delta$ (under the
identifications \eqref{eqn:identDh}).
\end{defn}

\begin{defn}\label{def:antipode}
An $\h$-Hopf algebroid is an $\h$-bialgebroid $A$ equipped with a
$\C$-linear map $S:A\to A$, the antipode, such that
    \begin{equation*}
    \begin{split}
&S(\mu_r(f)a)=S(a)\mu_l(f),\; S(a\mu_l(f))=\mu_r(f)S(a),\;
\mbox{for all }a\in A, f\in M_{\h^*},\\
 &m\circ (id\otimes S)\circ\Delta(a)= \mu_l(\varepsilon(a)1),\; \mbox{for all } a\in A,\\
 &m\circ(S\otimes id)\circ\Delta(a)=\mu_r(T_\alpha
 (\varepsilon(a)1)), \; \mbox{for all } a\in A_{\alpha\beta},
    \end{split}
    \end{equation*}
where $m:A\tilde\otimes A\to A$ denotes the multiplication and
$\varepsilon(a)1$ is the result of applying the difference
operator $\varepsilon(a)$ to the constant function $1\in
M_{\h^*}$.
\end{defn}

If there exists an antipode on an $\h$-bialgebroid, it is unique.
Furthermore, the antipode is anti-multiplicative,
anti-comultiplicative, unital, counital and interchanges the
moment maps $\mu_l$ and $\mu_r$, see \cite[Prop.
2.2]{KoelinkRosengren}.

\vspace{\baselineskip}

The FRST-construction provides many examples of $\h$-bialgebroids,
see \cite{EtingofVarchenko}, \cite{EtingofSchiffmann},
\cite{FelderVarchenko}, \cite{KoelinkRosengren}. We recall the
construction.

Let $\h$ and $M_{\h^*}$ be as before,
$V=\bigoplus_{\alpha\in\h^*}V_\alpha$ be a finite-dimensional
diagonalizable $\h$-module and $R:\h^*\times\C \to
\mathrm{End}_\h(V\otimes V)$ a meromorphic function that commutes
with the $\h$-action on $V\otimes V$. Let $\{e_x\}_{x\in X}$ be a
homogeneous basis of $V$, where $X$ is an index set. Write $R^{a
b}_{x y}(\lambda,z)$ for the matrix elements of $R$,
    \begin{eqnarray*}
R(\lambda,z)(e_a\otimes e_b)=\sum_{x,y\in X} R^{a b}_{x
y}(\lambda, z) e_x\otimes e_y,
    \end{eqnarray*}
and define $\omega:X\to \h^*$ by $e_x\in V_{\omega(x)}$. Let $A_R$
be the unital complex associative algebra generated by the
elements $\{L_{xy}(z)\}_{x,y\in X}$, with $z\in\C$, together with
two copies of $M_{\h^*}$, embedded as subalgebras. The elements of
these two copies will be denoted by $f(\lambda)$ and $f(\mu)$,
respectively. The defining relations of $A_R$ are
$f(\lambda)g(\mu) = g(\mu)f(\lambda)$,
    \begin{equation}
    \begin{split}
 f(\lambda)L_{x y}(z) = L_{x y}(z) f(\lambda+\omega(x)),\quad
 f(\mu)L_{x y}(z) =L_{x y}(z) f(\mu+\omega(y)),
    \end{split}
    \end{equation}
for all $f$, $g\in M_{\h^*}$, together with the $R L L$-relations
    \begin{eqnarray}\label{eqn:RLL}
\sum_{x,y\in X} R^{x y}_{a c}(\lambda, z_1/z_2)L_{x b}(z_1)L_{y
d}(z_2)= \sum_{x,y\in X} R_{x y}^{b d}(\mu, z_1/z_2)L_{c
y}(z_2)L_{a x}(z_1),
    \end{eqnarray}
for all $z_1$, $z_2\in\C$ and $a$, $b$, $c$, $d\in X$.

The bigrading on $A_R$ is defined by $L_{x y}(z)\in
A_{\omega(x),\omega(y)}$ and $f(\lambda)$, $f(\mu)\in A_{0,0}$.
The moment maps defined by
    \begin{eqnarray}
\mu_l(f)=f(\lambda),\qquad \mu_r(f)=f(\mu),
    \end{eqnarray}
make $A_R$ into a $\h$-algebra. The $\h$-invariance of $R$ ensures
that the bigrading is compatible with the $RLL$-relations
(\ref{eqn:RLL}).

Finally the co-unit and co-multiplication defined by
    \begin{gather}
\varepsilon(L_{a b}(z))=\delta_{a b} T_{-\omega(a)},\quad
\varepsilon(f(\lambda))=\varepsilon(f(\mu))=f, \\
 \Delta(L_{a b}(z))=\sum_{x\in X} L_{a x}(z)\otimes L_{x b}(z), \quad
\Delta(f(\lambda))=f(\lambda)\otimes 1, \quad
\Delta(f(\mu))=1\otimes f(\mu),
    \end{gather}
equip $A_R$ with the structure of an $\h$-bialgebroid, see
\cite{EtingofVarchenko}.

\subsection{Elliptic $U(2)$ quantum group}

We now give the results of the generalized FRST-construction when
applied to the $R$-matrix (\ref{eqn:ellR}), see
\cite{FelderVarchenko}. Let $0<q< 1$, $0<p< 1$. We assume that
$p$, $q$ are generic, it suffices to take $p$ and $q$
algebraically independent over $\Q$. We denote the corresponding
$\h$-bialgebroid by $\F_R(M(2))$, it is an elliptic analogue of
the algebra of polynomials on the space of complex $2\times
2$-matrices. The four $L$-generators will be denoted by
$\alpha(z)=L_{1,1}(z)$, $\beta(z)=L_{1,-1}(z)$,
$\gamma(z)=L_{-1,1}(z)$ and $\delta(z)=L_{-1,-1}(z)$.

\begin{defn}\label{def:ellqg}
The $\h$-algebroid $\F_R(M(2))$ is generated by $\alpha(z)$,
$\beta(z)$, $\gamma(z)$, $\delta(z)$, $z\in\C$, together with two
copies of $M_{\h^*}$ denoted by $f(\lambda)$, $f(\mu)$. The
defining relations are
    \begin{subequations} \label{eqn:ellRLL}
    \begin{align}
    \begin{split}\label{eqn:RLLabcd}
&\alpha(z_1)\alpha(z_2) = \alpha(z_2)\alpha(z_1),\quad
 \beta(z_1)\beta(z_2) = \beta(z_2)\beta(z_1),\\
&\gamma(z_1)\gamma(z_2) = \gamma(z_2)\gamma(z_1),\quad
 \delta(z_1)\delta(z_2) = \delta(z_2)\delta(z_1),
    \end{split}\\
    \begin{split}\label{eqn:RLLab}
 \alpha(z_1)\beta(z_2) = a(\mu,z_{12})\beta(z_2)\alpha(z_1)+
c(\mu,z_{12})\alpha(z_2)\beta(z_1),\\
 \beta(z_1)\alpha(z_2) = b(\mu,z_{12})\beta(z_2)\alpha(z_1)+
d(\mu,z_{12})\alpha(z_2)\beta(z_1),
    \end{split}\\
    \begin{split}\label{eqn:RLLac}
 a(\lambda,z_{12})\alpha(z_1)\gamma(z_2)+b(\lambda,z_{12})\gamma(z_1)\alpha(z_2)=
\gamma(z_2)\alpha(z_1),\\
 c(\lambda,z_{12})\alpha(z_1)\gamma(z_2)+d(\lambda,z_{12})\gamma(z_1)\alpha(z_2)=
\alpha(z_2)\gamma(z_1),
    \end{split}\\
    \begin{split}\label{eqn:RLLcd}
 \gamma(z_1)\delta(z_2) = a(\mu,z_{12})\delta(z_2)\gamma(z_1)+
c(\mu,z_{12})\gamma(z_2)\delta(z_1),\\
 \delta(z_1)\gamma(z_2) = b(\mu,z_{12})\delta(z_2)\gamma(z_1)+
d(\mu,z_{12})\gamma(z_2)\delta(z_1),
    \end{split}\\
    \begin{split}\label{eqn:RLLbd}
 a(\lambda,z_{12})\beta(z_1)\delta(z_2)+b(\lambda,z_{12})\delta(z_1)\beta(z_2)=
\delta(z_2)\beta(z_1),\\
 c(\lambda,z_{12})\beta(z_1)\delta(z_2)+d(\lambda,z_{12})\delta(z_1)\beta(z_2)=
\beta(z_2)\delta(z_1),
    \end{split}\\
    \begin{split}\label{eqn:RLLadbc}
 a(\lambda,z_{12})\alpha(z_1)\delta(z_2)+b(\lambda,z_{12})\gamma(z_1)\beta(z_2)
   = a(\mu,z_{12})\delta(z_2)\alpha(z_1)+c(\mu,z_{12})\gamma(z_2)\beta(z_1),\\
 c(\lambda,z_{12})\alpha(z_1)\delta(z_2)+d(\lambda,z_{12})\gamma(z_1)\beta(z_2)
   = a(\mu,z_{12})\beta(z_2)\gamma(z_1)+c(\mu,z_{12})\alpha(z_2)\delta(z_1),\\
 a(\lambda,z_{12})\beta(z_1)\gamma(z_2)+b(\lambda,z_{12})\delta(z_1)\alpha(z_2)
   = b(\mu,z_{12})\delta(z_2)\alpha(z_1)+d(\mu,z_{12})\gamma(z_2)\beta(z_1),\\
 c(\lambda,z_{12})\beta(z_1)\gamma(z_2)+d(\lambda,z_{12})\delta(z_1)\alpha(z_2)
  =  b(\mu,z_{12})\beta(z_2)\gamma(z_1)+d(\mu,z_{12})\alpha(z_2)\delta(z_1),
    \end{split}
    \end{align}
    \end{subequations}
with $z_{12}=z_1/z_2$, together with $f(\lambda)g(\mu) =
g(\mu)f(\lambda)$,
    \begin{equation}\label{eqn:commMh}
    \begin{split}
  f(\lambda)\alpha(z) =\alpha(z) f(\lambda+1),\quad
      &f(\mu)\alpha(z) =\alpha(z) f(\mu+1),\\
  f(\lambda)\beta(z) = \beta(z) f(\lambda+1) ,\quad
      &f(\mu)\beta(z) =\beta(z) f(\mu-1),\\
  f(\lambda)\gamma(z) = \gamma(z) f(\lambda-1) ,\quad
      &f(\mu)\gamma(z) =\gamma(z) f(\mu+1),\\
  f(\lambda)\delta(z) = \delta(z) f(\lambda-1),\quad
      &f(\mu)\delta(z) =\delta(z) f(\mu-1).
    \end{split}
    \end{equation}
The bigrading $\F_R(M(2))=\bigoplus_{m,n\in \Z, m+n\in 2\Z}
\F_{mn}$ is defined on the generators by
    \begin{equation*}
\alpha(z)\in\F_{1,1},\; \beta(z)\in\F_{1,-1},\;
\gamma(z)\in\F_{-1,1},\; \delta(z)\in\F_{-1,-1},\; f(\lambda),
f(\mu)\in\F_{0,0}.
    \end{equation*}
The comultiplication $\Delta: \F_R(M(2))\to \F_R(M(2))\tilde
\otimes \F_R(M(2))$ and counit $\varepsilon: \F_R(M(2))\to D_{\h}$
are algebra homomorphisms defined on the generators by
    \begin{equation}
    \begin{split}\label{eqn:comultqg}
\Delta\alpha(z)=\alpha(z)\otimes\alpha(z)+\beta(z)\otimes\gamma(z),\qquad
 &\Delta\beta(z)=\alpha(z)\otimes\beta(z)+\beta(z)\otimes\delta(z),\\
 \Delta\gamma(z)=\gamma(z)\otimes\alpha(z)+\delta(z)\otimes\gamma(z),\qquad
 &\Delta\delta(z)=\gamma(z)\otimes\beta(z)+\delta(z)\otimes\delta(z),\\
 \Delta f(\lambda)=f(\lambda)\otimes 1,\qquad &\Delta f(\mu)=1\otimes
 f(\mu),
    \end{split}
    \end{equation}
and
    \begin{equation*}
    \begin{split}
\varepsilon(\alpha(z))= T_{-1},\qquad
 \varepsilon(\beta(z))=\varepsilon(\gamma(z))=0, \qquad
 \varepsilon(\delta(z))= T_1, \qquad
 \varepsilon(f(\lambda))=\varepsilon(f(\mu))=f.
    \end{split}
    \end{equation*}
\end{defn}

\begin{remark}
Since $a(\lambda,q^2)=c(\lambda,q^2)$ and
$b(\lambda,q^2)=d(\lambda,q^2)$ we see that the $R$-matrix
(\ref{eqn:ellR}) is singular for $z=q^2$. Using (\ref{eqn:theta})
we compute
    \begin{equation*}
\det\begin{pmatrix}
    a(\lambda,z) & b(\lambda,z)\\
    c(\lambda,z) & d(\lambda,z)
\end{pmatrix}
=q^2\frac{\theta(zq^{-2})}{\theta(zq^2)}.
    \end{equation*}
We find that $z=q^2$, up to powers of $p$, is the only zero of the
determinant of $R$. In the case $z_{12}=q^2$ the right hand side
of the relations in (\ref{eqn:RLLab}) are multiples of each other,
so this also holds for the left hand sides giving
$b(\mu,q^2)\alpha(q^2z)\beta(z)=a(\mu,q^2)\beta(q^2z)\alpha(z)$.
Simplifying this identity and doing the same for
(\ref{eqn:RLLac})-(\ref{eqn:RLLbd}) we find
    \begin{subequations}\label{eqn:RLLs}
    \begin{align}
\label{eqn:RLLabs} \theta(q^{-2\mu})\, \alpha(q^2z)\, \beta(z) &=
q^2\theta(q^{-2(\mu+2)})\, \beta(q^2z)\, \alpha(z),\\
\label{eqn:RLLacs} \gamma(z)\, \alpha(q^2z) &= \alpha(z)\,
\gamma(q^2z),\\
\label{eqn:RLLcds} \theta(q^{-2\mu})\, \gamma(q^2z)\, \delta(z) &=
q^2\theta(q^{-2(\mu+2)})\, \delta(q^2z)\, \gamma(z),\\
\label{eqn:RLLbds} \delta(z)\, \beta(q^2z) &= \beta(z)\,
\delta(q^2z).
    \end{align}
    \end{subequations}
From the relations (\ref{eqn:RLLadbc}) we obtain three independent
relations
    \begin{equation}\label{eqn:RLLadbcs}
    \begin{split}
a(\lambda,q^2)\alpha(q^2z)\delta(z)+b(\lambda,q^2)\gamma(q^2z)\beta(z)&=
a(\mu,q^2)[\gamma(z)\beta(q^2z)+\delta(z)\alpha(q^2z)],\\
 a(\lambda,q^2)\beta(q^2z)\gamma(z)+b(\lambda,q^2)\delta(q^2z)\alpha(z)&=
 b(\mu,q^2)[\gamma(z)\beta(q^2z)+\delta(z)\alpha(q^2z)],\\
 \alpha(z)\delta(q^2z)+\beta(z)\gamma(q^2z)&=\gamma(z)\beta(q^2z)+\delta(z)\alpha(q^2z).
    \end{split}
    \end{equation}

From (\ref{eqn:ellRelem}) we see that $a(\lambda,z)$,
$b(\lambda,z)$, $c(\lambda,z)$ and $d(\lambda,z)$ have a simple
pole for $z=q^{-2}$. The residual relations of (\ref{eqn:RLL}) are
the relations obtained by multiplying by $z_{12}-q^{-2}$ and
taking the limit $z_{12}\to q^{-2}$, see \cite{FelderVarchenko}.
By convention, we interpret (\ref{eqn:RLL}) so that these are also
suppose to hold. The residual relations of (\ref{eqn:RLLab}),
respectively (\ref{eqn:RLLac}), (\ref{eqn:RLLcd}),
(\ref{eqn:RLLbd}) are linearly dependent and simplify to
(\ref{eqn:RLLs}). The residual relations of (\ref{eqn:RLLadbc})
reduce to three independent relations of which two can be derived
from (\ref{eqn:RLLadbcs}). The independent relation can be written
as
    \begin{equation}\label{eqn:RLLres}
    \begin{split}
b(\mu,q^{2})\gamma(q^2z)\beta(z)-a(\mu,q^{2})\delta(q^2z)\alpha(z)
&= a(\lambda,q^{2})[\gamma(z)\beta(q^2z)-\alpha(z)\delta(q^2z)].
    \end{split}
    \end{equation}

Note that $\Delta$, $\varepsilon$ preserve the commutation
relations (\ref{eqn:RLLs})-(\ref{eqn:RLLres}).
\end{remark}

The following lemma is \cite[Theorem $13$]{FelderVarchenko}; where
it is stated without proof.
\begin{lemma}\label{lem:det}
The element
    \begin{equation*}
    \begin{split}
\det(z)&=\frac{F(\mu)}{F(\lambda)}\left[ \alpha(z)\delta(q^2
z)-\gamma(z)\beta(q^2z)\right]\\
 &=\frac{F(\mu)}{F(\lambda)}\left[\delta(z)\alpha(q^2z)-\beta(z)\gamma(q^2z)
 \right]\\
 &=\frac{q^\mu}{q^\lambda}\left[\frac{\theta(q^{-2(\mu+2)})}{\theta(q^{-2(\lambda+2)})}\delta(q^2z)\alpha(z)
 - \frac{\theta(q^{-2\mu})}{q^2\theta(q^{-2(\lambda+2)})}\gamma(q^2z)\beta(z)\right]\\
 &=\frac{q^\mu}{q^\lambda}\left[\frac{\theta(q^{-2\mu})}{\theta(q^{-2\lambda})}\alpha(q^2z)\delta(z)-
  \frac{q^2\theta(q^{-2(\mu+2)})}{\theta(q^{-2\lambda})}\beta(q^2z)\gamma(z)\right],
    \end{split}
    \end{equation*}
where $F(\mu)= q^{\mu}\theta(q^{-2(\mu+1)})$, is a central element
of $\F_R(M(2))$. Moreover, $\Delta(\det(z))=\det(z)\otimes\det(z)$
and $\varepsilon(\det(z))=1$.
\end{lemma}

\begin{proof}
To see that the four expressions are equal we note that the first
equality is just the third equation of (\ref{eqn:RLLadbcs}), the
last equality follows from the equality obtained by eliminating
the right hand side in the first two equations of
(\ref{eqn:RLLadbcs}). The equality of the first and the third
expression is (\ref{eqn:RLLres}).

To prove that the element $\det(z)$ is central we have to show
that it commutes with every generator of the algebra, i.e.~
$\det(z)$ commutes with $\alpha(w)$, $\beta(w)$, $\gamma(w)$,
$\delta(w)$, $f(\mu)$ and $f(\lambda)$ for all $z$, $w\in\C$. We
write down the proof of $\beta(w)\det(z)=\det(z)\beta(w)$ in
detail, the other relations can be proved analogously. Using
(\ref{eqn:RLLabcd}), (\ref{eqn:RLLbd}) and (\ref{eqn:RLLab}) we
have
    \begin{equation*}
    \begin{split}
\beta(w)&\det(z) =F(\mu+1)/F(\lambda-1)
\left[\beta(w)\delta(z)\alpha(q^2z) - \beta(z)\beta(w)\gamma(q^2z)
\right]\\
 =& \,F(\mu+1)/F(\lambda-1)\left[c(\lambda,z/w) \beta(z)\delta(w)\alpha(q^2z) +
 d(\lambda,z/w)\delta(z)\beta(w)\alpha(q^2z)-\beta(z)\beta(w)\gamma(q^2z)\right]\\
 =& \,F(\mu+1)/F(\lambda-1)\left[c(\lambda,z/w) \beta(z)\delta(w)\alpha(q^2z)
 -\beta(z)\beta(w)\gamma(q^2z)\right.\\
 &\left.+ d(\lambda,z/w)\{b(\mu+1,w/(q^2z))\delta(z)\beta(q^2z)\alpha(w)
  + d(\mu+1,w/(q^2z))\delta(z)\alpha(q^2z)\beta(w)\}\right]\\
 =& \,F(\mu+1)/F(\lambda-1)\left[c(\lambda,z/w) \beta(z)\delta(w)\alpha(q^2z)
 -\beta(z)\beta(w)\gamma(q^2z)\right.\\
  &\left. + d(\lambda,z/w)\{b(\mu+1,w/(q^2z))\beta(z)\delta(q^2z)\alpha(w)
 +d(\mu+1,w/(q^2z))\delta(z)\alpha(q^2z)\beta(w))\} \right],
    \end{split}
    \end{equation*}
where we used (\ref{eqn:RLLbds}) in the last step. Since
    \begin{equation}\label{eqn:hulp}
-c(\lambda+1, z/w)= \frac{b(\lambda,w/(q^2z))}
{a(\lambda,w/(q^2z))}, \qquad
d(\lambda+1,z/w)=\frac{1}{a(\lambda,w/(q^2z))},
    \end{equation}
 we can simplify the last expression using the third relation of (\ref{eqn:RLLadbc}),
    \begin{equation*}
    \begin{split}
\beta(w)\det(z)=& \,F(\mu+1)/F(\lambda-1)\left[
-\frac{d(\mu+1,w/(q^2z))}{a(\lambda,w/(q^2z))}\beta(z)\gamma(q^2z)\beta(w)\right.\\
& \left.+ d(\mu+1,w/(q^2z))d(\lambda,z/w)
\delta(z)\alpha(q^2z)\beta(w) \right]\\
 =& \,\det(z)\beta(w),
    \end{split}
    \end{equation*}
where we use again the second relation of (\ref{eqn:hulp}) and
    \begin{equation*}
F(\mu+1)F(\lambda)d(\mu+1,w/(q^2z))d(\lambda,z/w)=F(\mu)F(\lambda-1),
    \end{equation*}
in the last step. From the definition of $\varepsilon$ it follows
that $\varepsilon (\det(z))=1$. Furthermore,
    \begin{equation*}
    \begin{split}
\Delta(\det(z))  =& \left(\frac{1}{F(\lambda)}\otimes
F(\mu)\right)\left( [\alpha(z)\delta(q^2z)-
\gamma(z)\beta(q^2z)]\otimes \alpha(z)\delta(q^2z)\right.\\
 &\qquad\left. +[\beta(z)\gamma(q^2z)-\delta(z)\alpha(q^2z)]\otimes
\gamma(z)\beta(q^2z) \right)\\
 =& \left(\frac{1}{F(\lambda)}\otimes
F(\mu)\right)\left(F(\lambda)/F(\mu)\det(z)\otimes
[\alpha(z)\delta(q^2z)-\gamma(z)\beta(q^2z)] \right) \\
 =&\, \det(z)\otimes\det(z),
    \end{split}
    \end{equation*}
where we use (\ref{eqn:RLLacs}) and (\ref{eqn:RLLbds}) in the
first equality.
\end{proof}

To $\F_R(M(2))$ we adjoin the central element $\det^{-1}(z)$
subject to the relation $\det(z)\det^{-1}(z)=1$. The
comultiplication and counit extend by
    \begin{equation*}
\Delta(\det^{-1}(z))=\det^{-1}(z)\otimes \det^{-1}(z), \qquad
\varepsilon(\det^{-1}(z))=1.
    \end{equation*}
It is easily checked that the resulting algebra, denoted by
$\F_R(GL(2,\C))$, is an $\h$-bialgebroid. Note that $\det(z)$ and
$\det^{-1}(z)$ have $(0,0)$-bigrading. In the dynamical
representations of Proposition \ref{prop:dynrepr} we consider
later, $\det(z)$ does not act as $id$, see Remark \ref{rem:det=1}.
Therefore we do not put $\det(z)=1$.

\begin{lemma}\label{lem:antipode}
The $\h$-bialgebroid $\F_R(GL(2,\C))$ is an $\h$-Hopf algebroid
with the antipode $S$ defined by $S(\det^{-1}(z))=\det(z)$,
    \begin{equation}
    \begin{split}\label{eqn:antipode}
S(\alpha(z))= \frac{F(\mu)}{F(\lambda)}
\det^{-1}(q^{-2}z)\delta(q^{-2}z), \quad
    & S(\beta(z)) = -\frac{F(\mu)}{F(\lambda)}\det^{-1}(q^{-2}z)\beta(q^{-2}z),\\
 S(\gamma(z)) = -\frac{F(\mu)}{F(\lambda)}\det^{-1}(q^{-2}z)\gamma(q^{-2}z), \quad
 & S(\delta(z)) = \frac{F(\mu)}{F(\lambda)}\det^{-1}(q^{-2}z)\alpha(q^{-2}z),\\
 S(f(\lambda))= f(\mu),\quad & S(f(\mu))= f(\lambda),
    \end{split}
    \end{equation}
on the generators and extended as an algebra antihomomorphism.
\end{lemma}

\begin{proof}
By Proposition 2.2 of \cite{KoelinkRosengren} we only have to
check that on the generators we have
    \begin{equation}\label{eqn:propS}
\begin{pmatrix}
S(\alpha(z))& S(\beta(z))\\ S(\gamma(z)) & S(\delta(z))
\end{pmatrix}
\begin{pmatrix}
\alpha(z)& \beta(z)\\ \gamma(z) & \delta(z)
\end{pmatrix}
=
\begin{pmatrix} 1& 0\\ 0 & 1
\end{pmatrix}
= \begin{pmatrix} \alpha(z)& \beta(z)\\ \gamma(z) & \delta(z)
\end{pmatrix}
\begin{pmatrix}
S(\alpha(z))& S(\beta(z))\\ S(\gamma(z)) & S(\delta(z))
\end{pmatrix},
    \end{equation}
and that the antipode preserves the defining relations of the
algebra. The proof is straightforward, using the $RLL$-relations
and Lemma \ref{lem:det}. Note that we need the residual relation
(\ref{eqn:RLLres}) for the second equality in (\ref{eqn:propS}).
\end{proof}

Next we give a $*$-structure to the obtained $\h$-Hopf algebroid.
Therefore we recall the definition of a $*$-structure on
$\h$-bialgebroids, see \cite{KoelinkRosengren}. Assuming
$\bar\;:\h\to\h$ is a conjugation, we put $\overline{f}(\lambda)=
\overline{f(\overline \lambda)}$, $f\in M_{\h^*}$. A $*$-operator
on an $\h$-bialgebroid $A$ is a $\C$-antilinear,
antimultiplicative involution on $A$ satisfying $\mu_l(\overline
f)=\mu_l(f)^*$, $\mu_r(\overline f)=\mu_r(f)^*$ and $(*\otimes
*)\circ \Delta=\Delta\circ *$, $\varepsilon\circ
*=*^{D_{\h}}\circ \varepsilon$ where $*^{D_{\h}}$ is defined by
$(f T_\alpha)^*= T_{-\overline\alpha}\overline f$. We use complex
conjugation on $\h \cong \C$.

\begin{lemma}\label{lem:*structure}
The $\h$-Hopf algebroid $\F_R(GL(2,\C))$ has a $*$-structure
defined on the generators by
$\det^{-1}(z)^*=\det^{-1}(q^{-2}/\overline z)$,
    \begin{equation}
    \begin{split}
\alpha(z)^*=\delta(1/\overline{z}),\quad
\beta(z)^*=-\gamma(1/\overline{z}),\quad&
 \gamma(z)^*=-\beta(1/\overline{z}),\quad
 \delta(z)^*=\alpha(1/\overline{z}).
    \end{split}
    \end{equation}
\end{lemma}

We call this $\h$-Hopf algebroid the elliptic $U(2)$ quantum group
and denote it by $\F_R(U(2))$.

\begin{proof}
We can easily check that this definition preserves the defining
relations of the algebra, that it is an involution and that we
have $(*\otimes*)\circ\Delta = \Delta \circ
*$ and $\varepsilon\circ
* = *^{D_\h}\circ \varepsilon$.
\end{proof}

\begin{remark}
Note that $S(\det(z))=\det^{-1}(z)$ and
$\det(z)^*=\det(q^{-2}/\overline z)$.
\end{remark}

\section{Corepresentations of the elliptic $U(2)$ quantum group}

Before discussing a special corepresentation of the elliptic
$U(2)$ quantum group, we recall the general definition of a
corepresentation of an $\h$-bialgebroid on an $\h$-space, see
\cite{KoelinkRosengren}.

\begin{defn}
An $\h$-space is a vector space over $M_{\h^*}$ which is also a
diagonalizable $\h$-module, $V=\bigoplus_{\alpha\in\h^*}V_\alpha$,
with $M_{\h^*}V_\alpha\subseteq V_\alpha$ for all $\alpha\in\h^*$.
A morphism of $\h$-spaces is an $\h$-invariant (i.e.~ grade
preserving) $M_{\h^*}$-linear map.
\end{defn}

We next define the tensor product of an $\h$-bialgebroid $A$ and
an $\h$-space $V$. Put  $A\tilde\otimes V=$ \linebreak
$\bigoplus_{\alpha,\beta\in\h^*}(A_{\alpha\beta}\otimes_{M_{\h^*}}V_\beta)$
where $\otimes_{M_{\h^*}}$ denotes the usual tensor product modulo
the relations $\mu_r^A(f) a\otimes v=a\otimes f v$. The grading
$A_{\alpha\beta}\otimes_{M_{\h^*}}V_\beta \subseteq
(A\tilde\otimes V)_\alpha$ and the extension of scalars
$f(a\otimes v)=\mu_l^A(f)a\otimes v$, $a\in A$, $v\in V$, $f\in
M_{\h^*}$, make $A\tilde\otimes V$ into an $\h$-space.

\begin{defn}
A (left) corepresentation of an $\h$-bialgebroid $A$ on an
$\h$-space $V$ is an $\h$-space morphism $\rho:V\to A\tilde\otimes
V$ such that
    \begin{equation}\label{eqn:propcorep}
(\Delta\otimes id)\circ \rho=(id\otimes\rho)\circ\rho, \qquad
(\varepsilon\otimes id)\circ \rho=id.
    \end{equation}
The first equality is in the sense of the natural isomorphism
$(A\tilde\otimes A)\tilde\otimes V\cong A\tilde\otimes
(A\tilde\otimes V)$ and in the second identity we use the
identification $V\simeq D_{\h}\tilde\otimes V$ defined by
$fT_{-\alpha}\otimes v\cong f v$, $f\in M_{\h^*}$, for all $v\in
V_\alpha$.
\end{defn}

Choose a homogeneous basis $\{v_k\}_k$ of $V$ over $M_{\h^*}$,
$v_k\in V_{\omega(k)}$, and introduce the corresponding matrix
elements of a corepresentation $\rho$ by $\rho(v_k)=\sum_j
t_{kj}\otimes v_j$. For these matrix elements we have from
(\ref{eqn:propcorep})
    \begin{equation*}
\Delta(t_{ij})=\sum_k t_{ik}\otimes t_{kj}, \qquad
\varepsilon(t_{ij})=\delta_{ij}T_{-\omega(i)},
    \end{equation*}
and Definition \ref{def:antipode} implies
    \begin{equation}\label{eqn:Sorth}
\delta_{kl}=\sum_j S(t_{kj})t_{jl} =\sum_j t_{kj}S(t_{jl}).
    \end{equation}

Our next objective is to construct explicit corepresentations of
$\F_R(U(2))$. Define, with the convention that the empty product
is $1$,
    \begin{equation}\label{eqn:vk}
v_k= v_k(z)= \gamma(z)\gamma(q^2z)\cdots\gamma(q^{2(N-k-1)})
\alpha(q^{2(N-k)})\cdots\alpha(q^{2(N-1)}z), \qquad
k\in\{0,1,\ldots,N\},
    \end{equation}
and put $V_{2k-N}=\mu_l(M_{\h^*})v_k$, $V=V^N=\bigoplus_{k=0}^N
V_{2k-N}$. Then $V$ is an $\h$-space with the multiplication by
$f\in M_{\h^*}$ given by the left moment map $\mu_l$. Note that
the grading on $V$ is compatible with Definition \ref{def:ellqg}.
We show that $\Delta: V^N \to \F_R(U(2))\tilde\otimes V^N$ making
$V^N$ a corepresentation of $\F_R(U(2))$, see Theorem
\ref{thm:corepmatrix}.

We start with the following preparatory lemma.

\begin{lemma}\label{lem:alphabetareverse}
In the $\h$-Hopf algebroid $\F_R(U(2))$ we have
    \begin{equation*}
    \begin{split}
\alpha(q^{2k}z)\beta(q^{2(l-1)}z)\cdots\beta(z)=\,&
\frac{\theta(q^{2(k-l+1)}, q^{2(\mu+l+1)})}{\theta(q^{2(k+1)},
q^{2(\mu+1)})} \beta(z)\cdots\beta(q^{2(l-1)}z)\alpha(q^{2k}z) \\
 &+ \frac{\theta(q^{2},q^{2(\mu+k+1)})}{\theta(q^{2(k+1)},
q^{2(\mu+1)})} \sum_{i=0}^{l-1}\beta(z)\cdots
\alpha(q^{2i}z)\cdots\beta(q^{2(l-1)}z) \beta(q^{2k}z),
    \end{split}
    \end{equation*}
for all $k\geq l\geq 1$.
\end{lemma}

\begin{proof}
For $l=1$ this is (\ref{eqn:RLLab}). In order to provide for the
induction step we interchange the order of the $\beta$'s by
(\ref{eqn:RLLabcd}), use the case $l=1$ and then the induction
hypothesis with $z\mapsto q^2 z$, $k\mapsto k-1$ finishes the
proof using (\ref{eqn:commMh}) and (\ref{eqn:ellRelem}).
\end{proof}

\begin{thm}\label{thm:corepmatrix}
In the $\h$-Hopf algebroid $\mathcal{F}_R(U(2))$, with $v_k(z)$
defined by \eqref{eqn:vk}, we have
    \begin{equation}
    \begin{split}
\Delta(v_k(z))= \sum_{j=0}^N t^N_{kj}(\mu,z)\otimes v_j(z),
    \end{split}
    \end{equation}
where the matrix-elements $t_{kj}^N(\mu,z)$ are given by
    \begin{equation*}
    \begin{split}
t^N_{kj}(\mu,z)=&\sum_{l=\max(0,k+j-N)}^{\min(k,j)} \ellbin{k}{l}
\ellbin{N-k}{j-l}\frac{(q^{2(\mu+N-k-2j+l+2)})_l}{(q^{2(\mu+N-2j+2)})_l}
 \frac{(q^{2(\mu+l-j+2)})_{j-l}}{(q^{2(\mu+N-2j-k+2l+2)})_{j-l}}\\
    &\quad\times\gamma(q^{2(N-k-1)}z)\cdots\gamma(q^{2(N-j-k+l)}z)
\delta(q^{2(N-j-k+l-1)}z)\cdots\delta(z)\\
    &\quad\times\alpha(q^{2(N-1)}z)\cdots\alpha(q^{2(N-l)}z)
\beta(q^{2(N-l-1)}z)\cdots\beta(q^{2(N-k)}z).
    \end{split}
    \end{equation*}
\end{thm}

\begin{proof}
We first deal with the cases $k=N$ and $k=0$, and get the general
result from the homomorphism property of the comultiplication
$\Delta$.
\begin{claim}
For all $k\in \Z_{\geq 0}$,
    \begin{equation}
    \begin{split}\label{eqn:claima}
\Delta(\alpha(z)\cdots\alpha(q^{2(k-1)}z))=&\sum_{l=0}^k
C_{kl}(\mu)\alpha(q^{2(k-1)}z)\cdots
\alpha(q^{2(k-l)}z)\beta(q^{2(k-l-1)}z)\cdots\beta(z)\\
    & \qquad\otimes \gamma(z)\cdots \gamma(q^{2(k-l-1)}z)\alpha(q^{2(k-l)}z)\cdots
\alpha(q^{2(k-1)}z),
    \end{split}
    \end{equation}
where the coefficients $C_{kl}(\mu)\in M_{\h^*}$ are given by
    \begin{equation*}
    \begin{split}
C_{kl}(\mu)=\ellbin{k}{l}\frac{(q^{2(\mu-l+2)})_l}{(q^{2(\mu+k-2l+2)})_l}.
    \end{split}
    \end{equation*}
\end{claim}
Note that $C_{k,0}=C_{k,k}=1$. We prove the claim by induction on
$k$. For $k=1$ this is just Definition \ref{def:ellqg} of the
comultiplication on $\alpha(z)$. Assume that the claim is true for
$k$. Then we obtain from (\ref{eqn:comultqg}) and repeated
application of (\ref{eqn:RLLacs}) that
    \begin{equation*}
    \begin{split}
\Delta(\alpha(z)\cdots\alpha(q^{2k}z))=&
\Delta(\alpha(z)\cdots\alpha(q^{2(k-1)}z))\Delta(\alpha(q^{2k})z)\\
 =&\sum_{l=0}^k C_{kl}(\mu)
\alpha(q^{2(k-1)}z)\cdots\alpha(q^{2(k-l)}z)\beta(q^{2(k-l-1)}z)\cdots\beta(z)\beta(q^{2k}z)\\
    &\qquad\otimes\gamma(z)\cdots\gamma(q^{2(k-l-1)}z)\alpha(q^{2(k-l)}z)\cdots
\alpha(q^{2(k-1)}z)\gamma(q^{2k}z) \\
    &+ \sum_{l=0}^k C_{kl}(\mu) \alpha(q^{2(k-1)}z)\cdots\alpha(q^{2(k-l)}z)
\beta(q^{2(k-l-1)}z)\cdots\beta(z)\alpha(q^{2k}z)\\
    &\qquad\otimes\gamma(z)\cdots\gamma(q^{2(k-l-1)}z)\alpha(q^{2(k-l)}z)\cdots\alpha(q^{2(k-1)}z)
\alpha(q^{2k}z) \\
 =& \sum_{l=0}^{k} C_{k,l}(\mu)
\alpha(q^{2(k-1)}z)\cdots\alpha(q^{2(k-l)}z)
\beta(q^{2(k-l-1)}z)\cdots\beta(z)\beta(q^{2k}z)\\
    &\qquad\otimes\gamma(z)\cdots\gamma(q^{2(k-l)}z)\alpha(q^{2(k-l+1)}z)\cdots
\alpha(q^{2k}z) \\
     &+ \sum_{l=1}^{k+1} C_{k,l-1}(\mu)
\alpha(q^{2(k-1)}z)\cdots\alpha(q^{2(k-l+1)}z)
\beta(q^{2(k-l)}z)\cdots\beta(z)\alpha(q^{2k}z)\\
    &\qquad\otimes\gamma(z)\cdots\gamma(q^{2(k-l)}z)\alpha(q^{2(k-l+1)}z)\cdots
\alpha(q^{2k}z).
    \end{split}
    \end{equation*}
For $l=0$ and $l=k+1$, we have $C_{k+1,0}(\mu)=C_{k,0}(\mu)=1$ and
$C_{k+1,k+1}(\mu)=C_{k,k}(\mu)=1$ respectively. So it remains to
prove that for $1\leq l\leq k$ we have
    \begin{equation}
    \begin{split}\label{eqn:hulpclaima}
C_{k+1,l}(\mu)\alpha(q^{2k}z)&\cdots\alpha(q^{2(k-l+1)}z)\beta(q^{2(k-l)}z)\cdots
\beta(z)\\
     =&\; C_{k,l}(\mu)\alpha(q^{2(k-1)}z)\cdots\alpha(q^{2(k-l)}z)
\beta(q^{2(k-l-1)}z)\cdots\beta(z)\beta(q^{2k}z)\\
     & + C_{k,l-1}(\mu)
\alpha(q^{2(k-1)}z)\cdots\alpha(q^{2(k-l+1)}z)\beta(q^{2(k-l)}z)\cdots\beta(z)\alpha(q^{2k}z).
    \end{split}
    \end{equation}
Using
    \begin{equation*}
C_{k,l-1}(\mu)=C_{k,l}(\mu)\frac{\theta(q^{2l}, q^{2(\mu+k-2l+3)},
q^{2(\mu+k-2l+2)})}{\theta(q^{2(k-l+1)}, q^{2(\mu+k-l+2)},
q^{2(\mu-l+2)})},
    \end{equation*}
and (\ref{eqn:commMh}) we obtain that the right hand side of
(\ref{eqn:hulpclaima}) equals
    \begin{equation*}
    \begin{split}
    C_{kl}(\mu)&\alpha(q^{2(k-1)}z)\cdots\alpha(q^{2(k-l+1)}z)\\
      &\times\left[\alpha(q^{2(k-l)}z)\beta(q^{2(k-l-1)}z)\cdots\beta(z)\beta(q^{2k}z)\right.\\
     &\quad +\left. \frac{\theta(q^{2l}, q^{2(\mu+k-l+2)},
q^{2(\mu+k-l+1)})}{\theta(q^{2(k-l+1)}, q^{2(\mu+k+1)},
q^{2(\mu+1)})}\beta(q^{2(k-l)}z)\cdots\beta(z)\alpha(q^{2k}z)\right].
    \end{split}
    \end{equation*}
By Lemma \ref{lem:alphabetareverse}, with $(k,l)$ replaced by
$(k-l,k-l)$ the term in square brackets equals
    \begin{equation*}
    \begin{split}
\frac{\theta(q^2,q^{2(\mu+k-l+1)})}{\theta(q^{2(k-l+1)},q^{2(\mu+1)})}&
\sum_{n=0}^{k-l-1}\beta(z)\cdots\alpha(q^{2n}z)\cdots\beta(q^{2(k-l)}z)\beta(q^{2k}z)\\
    &+ \frac{\theta(q^2,q^{2(\mu+k-l+1)})}{\theta(q^{2(k-l+1)},q^{2(\mu+1)})}
\beta(z)\cdots\beta(q^{2(k-l-1)}z)\alpha(q^{2(k-l)}z)\beta(q^{2k}z)\\
 &\quad+ \frac{\theta(q^{2l}, q^{2(\mu+k-l+2)},
q^{2(\mu+k-l+1)})}{\theta(q^{2(k-l+1)}, q^{2(\mu+k+1)},
q^{2(\mu+1)})}\beta(z)\cdots\beta(q^{2(k-l)}z)\alpha(q^{2k}z),\\ =
\frac{\theta(q^{2(k+1)},
q^{2(\mu+k-l+1)})}{\theta(q^{2(k-l+1)},q^{2(\mu+k+1)})}&
\left[\frac{\theta(q^{2l}, q^{2(\mu+k-l+2)})}{\theta(q^{2(k+1)},
q^{2(\mu+1)})}
\beta(z)\cdots\beta(q^{2(k-l)z})\alpha(q^{2k}z)\right.\\
 & \quad+ \left.\frac{\theta(q^{2}, q^{2(\mu+k+1)})}{\theta(q^{2(k+1)},q^{2(\mu+1)})}\sum_{n=0}^{k-l}
\beta(z)\cdots\alpha(q^{2n}z)\cdots\beta(q^{2(k-l)}z)\beta(q^{2k}z)\right]\\
 = \frac{\theta(q^{2(k+1)},
q^{2(\mu+k-l+1)})}{\theta(q^{2(k-l+1)},q^{2(\mu+k+1)})}&
\alpha(q^{2k}z)\beta(q^{2(k-l)}z)\cdots\beta(z),
    \end{split}
    \end{equation*}
where we use Lemma \ref{lem:alphabetareverse} with $(k,l)$
replaced by $(k,k-l+1)$ in the last step. Using
(\ref{eqn:RLLabcd}) and (\ref{eqn:commMh}) we see that the right
hand side of (\ref{eqn:hulpclaima}) equals the left hand side
using
    \begin{equation*}
C_{k+1,l}(\mu)=
C_{k,l}(\mu)\frac{\theta(q^{2(k+1)},q^{2(\mu+k-2l+2)})}{\theta(q^{2(k-l+1)},
q^{2(\mu+k-l+2)})}.
    \end{equation*}
This proves the claim.

Since the ($\alpha$,$\beta$)- and the
($\gamma$,$\delta$)-commutation relations are similar by
(\ref{eqn:RLLab}), (\ref{eqn:RLLcd}), (\ref{eqn:RLLabs}),
(\ref{eqn:RLLcds}), we analogously have
    \begin{equation}
    \begin{split}\label{eqn:claimc}
\Delta(\gamma(z)\cdots\gamma(q^{2(k-1)}z))=&\sum_{l=0}^k
C_{kl}(\mu)\gamma(q^{2(k-1)}z)\cdots
\gamma(q^{2(k-l)}z)\delta(q^{2(k-l-1)}z)\cdots\delta(z)\\
    &\qquad \otimes \gamma(z)\cdots
\gamma(q^{2(k-l-1)}z)\alpha(q^{2(k-l)}z)\cdots
\alpha(q^{2(k-1)}z).
    \end{split}
    \end{equation}

Using (\ref{eqn:claima}), (\ref{eqn:claimc}) and that the
comultiplication $\Delta$ is a morphism we find
    \begin{equation*}
    \begin{split}
\Delta(\gamma(z)&\cdots\gamma(q^{2(N-k-1)}z)\alpha(q^{2(N-k)}z)\cdots\alpha(q^{2(N-1)}z))\\
    =&\;\Delta(\gamma(z)\cdots\gamma(q^{2(N-k-1)}z))
\Delta(\alpha(q^{2(N-k)}z)\cdots\alpha(q^{2(N-1)}z))\\
    =&\;\sum_{l=0}^k\sum_{m=0}^{N-k}
C_{N-k,m}(\mu)C_{k,l}(\mu-2m+N-k)
    \gamma(q^{2(N-k-1)}z)\cdots\gamma(q^{2(N-k-m)}z)\\
     &\times \delta(q^{2(N-k-m-1)}z)\cdots\delta(z)\alpha(q^{2(N-1)}z)\cdots
\alpha(q^{2(N-l)}z)\beta(q^{2(N-l-1)}z)\cdots\beta(q^{2(N-k)}z)\\
    &\qquad\otimes
\gamma(z)\cdots\gamma(q^{2(N-m-l-1)}z)\alpha(q^{2(N-m-l)}z)\cdots\alpha(q^{2(N-1)}z),
    \end{split}
    \end{equation*}
where we use (\ref{eqn:commMh}), (\ref{eqn:RLLacs}). Substituting
$m=j-l$ gives
    \begin{equation*}
    \begin{split}\label{eqn:tkjhulp}
t^N_{kj}(\mu,z)=&\sum_{l=\max(0,j+k-N)}^{\min(k,j)}
C_{N-k,j-l}(\mu)C_{k,l}(\mu+N-2j+2l-k)\\
    &\times\gamma(q^{2(N-k-1)}z)\cdots
\gamma(q^{2(N-k-j+l)}z)\delta(q^{2(N-k-j+l-1)}z)\cdots\delta(z)\\
    & \times\alpha(q^{2(N-1)}z)\cdots
\alpha(q^{2(N-l)}z)\beta(q^{2(N-l-1)}z)\cdots\beta(q^{2(N-k)}z),
    \end{split}
    \end{equation*}
which proves the theorem.
\end{proof}

In the next proposition we prove that this corepresentation is
unitary in a certain sense. Note that this property is an
extension of unitarizability of a corepresentation introduced in
\cite{KoelinkRosengren}.

\begin{prop}\label{prop:unitarity}
The matrix elements $t_{kj}^N(\mu,z)$ of the corepresentation in
Theorem \ref{thm:corepmatrix} satisfy
    \begin{equation*}
\Gamma_k(\mu)
S(t_{kj}^N(\mu,z))^*=\Gamma_j(\lambda)t_{jk}^N(\mu,q^{-2(N-2)}/\overline
z)\prod_{i=0}^{N-1}\det^{-1}(q^{-2i}/\overline z),
    \end{equation*}
with
    \begin{equation*}
\Gamma_k(\mu)=
\ellbin{N}{k}\frac{(q^{2(\mu-k+2)})_k}{(q^{2(\mu+N-2k+2)})_k}
    \prod_{i=0}^{N-k-1}\frac{q^{-(\mu+N-2k-i)}}{\theta(q^{-2(\mu+N-2k-i+1)})}
    \prod_{i=0}^{k-1}\frac{q^{-(\mu-k+i)}}{\theta(q^{-2(\mu-k+i+1)})}.
    \end{equation*}
\end{prop}

\begin{proof}
To simplify the formulas in the proof we denote
$D=\prod_{i=0}^{N-1} \det^{-1}(q^{-2i}/\overline z)$,
$G_{Nk}(\mu)={\ellbin{N}{k}}\frac{(q^{2(\mu-k+2)})_k}{(q^{2(\mu-N-2k+2)})_k}$
and $F_{k}(\mu)=\prod_{i=0}^{k-1}F(\mu+i)$ where $F$ is defined in
Lemma \ref{lem:det}.

From Theorem \ref{thm:corepmatrix} we see that the matrix elements
$t_{kj}^N(\mu,z)$ for $k$ or $j$ equal to $0$ or $N$ consist of a
single term. Using Lemmas \ref{lem:antipode}, \ref{lem:*structure}
and the relations in Definition \ref{def:ellqg} proves the
proposition in case $j=N$
    \begin{equation*}
    \begin{split}
[S(t_{kN}^N(\mu,z))]^*= & D\,
\frac{F_{N-k}(\mu-k+1)}{F_{N-k}(\lambda-N)}
\frac{F_k(\lambda-k)}{F_k(\mu-k)} \\
     &\times \alpha(q^2/\overline
z)\cdots\alpha(q^{-2(k-2)}/\overline z)
\beta(q^{-2(k-1)}/\overline z)\cdots\beta(q^{-2(N-2)}/\overline
z)\\
 =& D\,\frac{F_{N-k}(\mu-k+1)}{F_{N-k}(\lambda-N)}
\frac{F_k(\lambda-k)}{F_k(\mu-k)} G_{Nk}(\mu)^{-1} t_{Nk}^N(\mu,
q^{-2(N-2)}/\overline z).
    \end{split}
    \end{equation*}
From $\Delta(t_{kN}^N(\mu,z))=\sum_{j=0}^N t_{kj}^N(\mu,z) \otimes
t_{jN}^N(\mu,z)$ and $\sigma\circ((*\circ S)\otimes (*\circ
S))\circ \Delta=\Delta\circ(*\circ S)$ we obtain
    \begin{equation*}
\sum_{j=0}^N S(t_{jN}^N(\mu,z))^*\otimes S(t_{kj}^N(\mu,z))^*=
\Delta(S(t_{kN}^N(\mu,z))^*).
    \end{equation*}
This relation gives
    \begin{equation*}
    \begin{split}
\sum_{j=0}^N  &F_{N-j}(\mu-j+1) F_j(\mu-j) G_{Nj}(\mu)^{-1}
t_{Nj}^N(\mu,q^{-2(N-2)}/\overline z)\otimes
S(t_{kj}^N(\mu,z))^*\\
  &= [1\otimes D \, F_{N-k}(\mu-k+1)
F_k(\mu-k) G_{Nk}(\mu)^{-1}]\sum_{j=0}^N
t_{Nj}^N(\mu,q^{-2(N-2)}/\overline z) \otimes
t_{jk}^N(\mu,q^{-2(N-2)}/\overline z).
    \end{split}
    \end{equation*}
Since $\{t_{Nj}^N(\mu,z)\}_{j=0}^N$ are linearly independent (this
follows easily from Proposition \ref{prop:dynrepr} and Lemma
\ref{lemma:pit}), the identity holds termwise. So
(\ref{eqn:tensor}) proves the proposition.
\end{proof}

\section{Discrete bi-orthogonality for elliptic hypergeometric series}

Using Proposition \ref{prop:unitarity} we can reformulate the
orthogonality relations (\ref{eqn:Sorth}) for the matrix elements
as
    \begin{subequations}\label{eqn:unorth}
    \begin{align}
    \label{eqn:unorth1}
\delta_{kl} = &\sum_{j=0}^N
(t_{jl}^N(\mu,z))^*\frac{\Gamma_j(\lambda)}{\Gamma_k(\mu)}
t_{jk}(\mu,q^{-2(N-2)}/\overline
z)\prod_{i=0}^{N-1}\det^{-1}(q^{-2i}/\overline z)\\
    \label{eqn:unorth2}
    =&\sum_{j=0}^N\frac{\Gamma_l(\lambda)}{\Gamma_j(\mu)}t^N_{lj}(\mu,q^{-2(N-2)}/\overline
z) (t_{kj}^N(\mu,z))^*\prod_{i=0}^{N-1}\det^{-1}(q^{-2i}/\overline
z).
    \end{align}\end{subequations}

To obtain commutative versions of (\ref{eqn:unorth}), we need to
represent the algebra $\F_R(U(2))$ explicitly. For this we need
the notion of a dynamical representation of an $\h$-algebra, see
\cite{EtingofSchiffmann}, \cite{EtingofVarchenko},
\cite{FelderVarchenko}, \cite{KoelinkRosengren}.

Let $V=\bigoplus_{\alpha\in\h^*} V_{\alpha}$ be an $\h$-space and
let $(D_{\h,V})_{\alpha\beta}$ be the space of $\C$-linear
operators $U$ on $V$ such that $U(gv)=T_{-\beta}(g)U(v)$ and
$U(V_\gamma)\subseteq V_{\gamma+\beta-\alpha}$ for all $g\in
M_{\h^*}$, $v\in V_\beta$, $\gamma\in\h^*$. Then the space
$D_{\h,V}=\bigoplus_{\alpha,\beta\in\h^*}(D_{\h,V})_{\alpha,\beta}$
is an $\h$-algebra with the moment maps $\mu_l$,
$\mu_r:M_{\h^*}\to (D_{\h,V})_{00}$ given by
$\mu_l(f)(v)=T_{-\alpha}(f)(v)$ and $\mu_r(f)(v)=fv$ for all $v\in
V_{\alpha}$.

\begin{defn}
A dynamical representation of an $\h$-algebra $A$ on an $\h$-space
$V$ is an $\h$-algebra homomorphism $A \to D_{\h,V}$.
\end{defn}

\begin{prop}\label{prop:dynrepr}
(see \cite{FelderVarchenko}) Let $\omega\in\C$ be arbitrary and
$\mathcal{H}^\omega$ be the $\h$-space with basis
$\{e_k\}_{k=0}^\infty$ and weight decomposition
$\mathcal{H}^\omega=\bigoplus_{k=0}^\infty
\mathcal{H}^\omega_{\omega-2k}$,
$\mathcal{H}^\omega_{\omega-2k}=M_{\h^*}e_k$. Then there exists a
dynamical representation $\pi^{\omega}:\F_R(M(2))\to
D_{\h,\mathcal{H}^\omega}$, defined on the generators by
    \begin{equation}
    \begin{split}
&\pi^\omega(\alpha(z))(g e_k)=A_k(\lambda,z)T_{-1}g e_k, \qquad
\pi^\omega(\beta(z))(g e_k)=B_k(\lambda,z)T_{1}g e_{k+1}\\
 &\pi^\omega(\gamma(z))(g e_k)=C_k(\lambda,z)T_{-1}g
 e_{k-1},\qquad
 \pi^\omega(\delta(z))(g e_k)=D_k(\lambda,z)T_{1}g e_k\\
 &\pi^\omega(\mu_r(f))(g e_k)=f(\lambda)g e_k, \qquad
 \pi^\omega(\mu_l(f))(g e_k)=f(\lambda-\omega+2k)g e_k,
    \end{split}
    \end{equation}
where $g\in M_{\h^*}$ and
    \begin{equation*}
    \begin{split}
A_k(\lambda,z) &= q^{2k} \frac{\theta(q^{-2(\lambda+1)-2k})\theta(
zq^{\omega-2k+1})} {\theta(q^{-2(\lambda+1)})\theta(
zq^{\omega+1})},\\
 B_k(\lambda,z) &= q^k\frac{\theta(q^2)\theta( zq^{-2(\lambda+1)+\omega-2k-1})}
  {\theta(q^{-2(\lambda+1)})\theta( zq^{\omega+1})},\\
 C_k(\lambda,z) &= q^{-(k-1)}\frac{\theta(q^{2k})\theta(q^{2(\omega-k+1)})
  \theta( zq^{2(\lambda+1)-\omega+2k-1})}
  {\theta(q^2)\theta(q^{2(\lambda+1)})\theta( zq^{\omega+1})}, \; C_0(\lambda,z)=0, \\
 D_k(\lambda,z) &= \frac{\theta(q^{-2(\lambda+1-\omega+k)})\theta( zq^{-\omega+2k+1})}
  {\theta(q^{-2(\lambda+1)})\theta( zq^{\omega+1})}.
    \end{split}
    \end{equation*}
\end{prop}

\begin{proof}
The dynamical representation preserves the defining relations
(\ref{eqn:ellRLL}), (\ref{eqn:commMh}) of the algebra as can be
checked by use of (\ref{eqn:theta}).
\end{proof}

\begin{remark}\label{rem:det=1}
Using the addition formula (\ref{eqn:theta}) we obtain
\begin{equation*}
\pi^\omega(\det(z))= q^\omega\frac{\theta( zq^{1-\omega})}{\theta(
zq^{1+\omega})}id,
\end{equation*}
so $\det(z)$ acts as a scalar. Note that this scalar is $1$ if
$\omega =0$.
\end{remark}

The action of a matrix element in the dynamical representation can
be calculated in terms of elliptic hypergeometric series.

\begin{lemma}\label{lemma:pit}
For the dynamical representation of Proposition \ref{prop:dynrepr}
we have
    \begin{equation*}
\pi^\omega(t^N_{kj}(\mu,z))(ge_m)=\tau_{kjm}^{N\omega}(\lambda,z)
(T_{N-2j}g)e_{m+k-j},
    \end{equation*}
where $\tau^{N\omega}_{kjm}(\lambda,z)$ is given by
    \begin{equation*}
    \begin{split}
\tau^{N\omega}_{kjm}(\lambda,z)=&(-1)^{N-k}\theta(q^2)^{k-j}
q^{\frac{3}{2}k(k-1)+N(N+1)+\frac{5}{2}j(j+1)+2N(\lambda-k-2j)+m(k-j)+3jk-2k\lambda}\\
&\times\frac{(q^{-2(\lambda+1)}, q^{2(m+k-j+1)},
q^{2(N-k-j+1)},q^{2(\omega-m-k+1)},
zq^{2(\lambda+N-2j+m+2)-\omega-1})_j} {(q^2,
q^{2(\lambda+N-k-2j+2)}, q^{2(\lambda-j+2)})_j}\\
 &\times \frac{( zq^{-2(\lambda-2j+m+k)+\omega-1})_k}{(q^{-2(\lambda+N-2j)})_k}
 \frac{(q^{-2(\lambda+N-2j-\omega+m)}, zq^{2(m+k)-\omega+1})_{N-k-j}}{(q^{2(\lambda-j+1)})_{N-k-j}}
\frac{1}{( zq^{\omega+1})_N}\\
 &\times
{}_{10}\omega_9[q^{2(\lambda+N-2j-k+1)}; q^{-2k}, q^{-2j},
q^{2(\lambda-j+1)}, q^{2(\lambda+N-2j-\omega+m+1)},\\
 &\qquad q^{2(\lambda+N+2+m-2j)},  zq^{2(N-m-k)+\omega+1},
  z^{-1}q^{-2(m+k-1)+\omega-1}].
    \end{split}
    \end{equation*}
\end{lemma}

\begin{proof}
From Proposition \ref{prop:dynrepr} and Theorem
\ref{thm:corepmatrix} it follows
    \begin{equation}\label{eqn:tauABCD}
    \begin{split}
\pi^\omega(t^{N}_{kj}(\mu,z))(ge_m)=&\sum_{l=\max(0,k+j-N)}^{\min(k,j)}
\ellbin{k}{l}\ellbin{N-k}{j-l}
\frac{(q^{2(\lambda+N-k-2j+l+2)})_l}{(q^{2(\lambda+N-2j+2)})_l}
 \frac{(q^{2(\lambda+l-j+2)})_{j-l}}{(q^{2(\lambda+N-2j-k+2l+2)})_{j-l}}\\
 &\times \prod_{n=0}^{j-l-1}C_{m+k-l-n}(\lambda-j+l+1+n,
 q^{2(n+N-k-j+l)}z)\\
  & \times\prod_{n=0}^{N-k-j+l-1}D_{m+k-l}(\lambda+N-k-2j+2l-1-n,
 q^{2n}z)\\
  & \times\prod_{n=0}^{j-l-1}A_{m+k-l}(\lambda+n+N-k-2j+l+1,
  q^{2(N-l+n)}z)\\
 & \times\prod_{n=0}^{k-l-1}B_{m+n}(\lambda+N-2j-1-n, q^{2(n+N-k)}z)
 \,(T_{N-2j}g)e_{m+k-j}.
    \end{split}
    \end{equation}
This gives the required form of the lemma, and it remains to show
that we can identify $\tau_{kjm}^{N\omega}(\lambda,z)$ with an
elliptic hypergeometric series. From the explicit expressions of
Proposition \ref{prop:dynrepr} we see that we can rewrite the four
products in terms of elliptic factorials
    \begin{equation}\label{eqn:ABCD}
    \begin{split}
& \prod_{n=0}^{j-l-1}A_{m+k-l}(\lambda+n+N-k-2j+l+1,
  q^{2(N-l+n)}z)= (-1)^l z^{-l}q^{2l(l-N)-l(\omega+l)}\\
   &\qquad\times \frac{(q^{2(\lambda+N-2j+m+2)}, zq^{2(N-m-k)+\omega+1}, q^{2(\lambda+N-2j-k+2)})_l}
   {(q^{-2(N-1)-\omega-1}/z)_l (q^{2(\lambda+N-2j-k+2)})_{2l}},\\
& \prod_{n=0}^{k-l-1}B_{m+n}(\lambda+N-2j-1-n, q^{2(n+N-k)}z)=
(-1)^l q^{2l(m+k-l)+l(l+1)+\frac{1}{2}(k-l)(2m+k-l-1)}\\
  &\qquad\times \theta(q^2)^k\frac{(q^{2(\lambda+N-2j-k+1)}, q^{-2(N-1)-\omega-1}/z)_l
  (zq^{-2(\lambda+k+m-2j)+\omega-1})_k}
  {(q^{2(\lambda+m-2j+1)-\omega+1})_l (q^{-2(\lambda+N-2j)}, zq^{2(N-k)+\omega+1})_k},\\
&  \prod_{n=0}^{j-l-1}C_{m+k-l-n}(\lambda-j+l+1+n,
q^{2(n+N-k-j+l)}z)= (-1)^l \theta(q^2)^{-j}\\
    &\qquad\times q^{j(l+1-m-k)-l(m+k+2)+ \frac{1}{2}(j-l)(j-l-1)}
 \frac{(q^{2(\lambda-j+2)}, zq^{2(N-k-j)+\omega+1})_l}
  {(q^{-2(m+k)}, q^{2(\omega-m-k+1)}, zq^{2(\lambda+N-2j+m+2)-\omega-1})_l}\\
  &\qquad\times\frac{(q^{2(\omega-m-k+1)}, zq^{2(\lambda+N-2j+M+2)-\omega-1}, q^{2(m+k-j+1)})_j}
  {(q^{2(\lambda-j+2)},zq^{2(N-k-j)+\omega+1})_j},\\
   & \prod_{n=0}^{N-k-j+l-1}D_{m+k-l}(\lambda+N-k-2j+2l-1-n,
 q^{2n}z)= (-1)^{N-j-k+l}z^l \\
  &\qquad \times q^{(N-j-k)(2\lambda+N-3j-k+2l+2)+l(\omega+l)}
   \frac{(q^{-2(\lambda+N-2j-\omega+m)}, zq^{2(m+k)-\omega+1})_{N-k-j}}
   {(q^{2(\lambda-j+1)},zq^{\omega+1})_{N-k-j}}\\
   &\qquad \times \frac{(q^{2(\lambda+N-2j-\omega+m+1)},
  q^{2(\lambda-j+1)}, q^{-2(m+k-1)+\omega-1}/z)_l}
  {(q^{2(\lambda+N-2j-k+1)})_{2l}(zq^{2(N-k-j)+\omega+1})_l },
    \end{split}
    \end{equation}
where we use elementary transformation formulas for the elliptic
factorials including
     \begin{equation*}
    \begin{split}\label{eqn:transell}
  (a q^{-4l})_l=(-1)^l(aq^{-4l})^lq^{l(l-1)}\frac{(q^2/a)_{2l}}{(q^2/a)_l}.
    \end{split}
    \end{equation*}
Furthermore for the elliptic binomials and the other factor in
(\ref{eqn:tauABCD}) we have
\begin{equation}\label{eqn:ellbinfact}
\begin{split}
&\ellbin{k}{l}= (-1)^l q^{2l(k-l+1)+l(l-1)}
\frac{(q^{-2k})_l}{(q^2)_l},\\
 &\ellbin{N-k}{j-l}= (-1)^l q^{2l(j-l+1)+l(l-1)}\frac{(q^{-2j})_l(q^{2(N-k-j+1)})_j}
 {(q^{2(N-k-j+1)})_l(q^2)_j},\\
 &\frac{(q^{2(\lambda+N-k-2j+l+2)})_l}{(q^{2(\lambda+N-2j+2)})_l}
 \frac{(q^{2(\lambda+l-j+2)})_{j-l}}{(q^{2(\lambda+N-2j-k+2l+2)})_{j-l}}= (-1)^{j}q^{2j(\lambda-j+2)+j(j-1)}\\
 &\qquad \times \frac{[(q^{2(\lambda+N-k-2j+2)})_{2l}]^2(q^{-2(\lambda+1)})_j}
 {(q^{2(\lambda+N-k-2j+2)},q^{2(\lambda+N-2j+2)},
 q^{2(\lambda-j+2)},q^{2(\lambda+N-k-j+2)})_l(q^{2(\lambda+N-2j-k+2)})_j}.
\end{split}
\end{equation}
Then, substituting (\ref{eqn:ABCD}) and (\ref{eqn:ellbinfact})
into (\ref{eqn:tauABCD}) gives
    \begin{equation*}
    \begin{split}
\tau^{N\omega}_{kjm}(\lambda,z)=
&(-1)^{N-k}\theta(q^2)^{k-j}q^{\frac{3}{2}k(k-1)+N(N+1)+\frac{5}{2}j(j+1)+2N(\lambda-k-2j)+m(k-j)+3jk-2k\lambda}\\
    &\times\frac{(q^{-2(\lambda+1)}, q^{2(m+k-j+1)}, q^{2(N-k-j+1)},
q^{2(\omega-m-k+1)}, zq^{2(\lambda+N-2j+m+2)-\omega-1})_j} {(q^2,
q^{2(\lambda+N-k-2j+2)}, q^{2(\lambda-j+2)},
 zq^{2(N-k-j)+\omega+1})_j}\\
 &\times \frac{( zq^{-2(\lambda-2j+m+k)+\omega-1})_k}{(q^{-2(\lambda+N-2j)},
 zq^{2(N-k)+\omega+1})_k} \frac{(q^{-2(\lambda+N-2j-\omega+m)}, zq^{2(m+k)-\omega+1})_{N-k-j}}{(q^{2(\lambda-j+1)},
 zq^{\omega+1})_{N-k-j}}\\
 &\times\sum_{l=\max(0,k+j-N)}^{\min(k,j)} q^{2l}\frac{\theta(q^{2(\lambda+N-k-2j+1+2l)})}{\theta(q^{2(\lambda+N-k-2j+1)})}
\frac{(q^{2(\lambda+N-2j-k+1)})_l}{(q^{2})_l} \frac{(q^{-2k})_l}{(
q^{2(\lambda+N-2j+2)})_l}\\
 &\times \frac{( q^{-2j})_l}{( q^{2(\lambda+N-k-j+2)})_l}
\frac{( q^{2(\lambda-j+1)})_l} {(q^{2(N-k-j+1)})_l}
\frac{(q^{2(\lambda+N-2j-\omega+m+1)})_l}
{(q^{2(\omega-k+m+1)})_l}
\frac{(q^{2(\lambda+N+2+m-2j)})_l}{(q^{-2(m+k)})_l}\\ &\times
\frac{(  zq^{2(N-k-m)+\omega+1})_l}{(
   z^{-1}q^{2(\lambda-2j+m+1)-\omega+1})_l}
 \frac{( q^{-2(m+k-1)+\omega-1}/z)_l}{(
 zq^{2(\lambda+N-2j+m+2)-\omega-1})_l}.
    \end{split}
    \end{equation*}
Note that if $k+j\geq N$ one of the factors in the denominator of
the sum, $(q^{2(N-k-j+1)})_l$, equals zero. However this pole is
cancelled by $(q^{2(N-k-j+1)})_j$.

\end{proof}

Analogously we can compute
\begin{lemma}\label{lemma:pit*}
For the dynamical representations of Proposition
\ref{prop:dynrepr} we have
    \begin{equation*}
\pi^\omega((t^N_{kj}(\mu,z))^*)(ge_m)=\tilde\tau^{N\omega}_{kjm}(\lambda,z)(T_{-N+2j}g)e_{m+j-k},
    \end{equation*}
where $\tilde\tau^{N\omega}_{kjm}(\lambda,z)$ is given by
    \begin{equation*}
    \begin{split}
\tilde\tau^{N\omega}_{kjm}(\lambda,z)=&
q^{2j^2+\frac{1}{2}k(k-1)-\frac{1}{2}j(j-1)+k(1-m-j)+2m(N-j-k)+mj+2(N-k)(\lambda-k+1)-(N-k-j)(N-k-j-1)}\\
 &\times (-1)^{N-j} \theta(q^2)^{j-k}
\frac{(q^{2(N-k-j+1)},q^{-2(\lambda-N+2j+1)},
q^{-2(\lambda+2j-k+m)+\omega+1}/\overline z)_j}
{(q^2)_j(q^{2(\lambda-k+2)},q^{-2(\lambda+2j-N)},
q^{-2(N-k-1)+\omega+1}/\overline z )_j}\\
 &\times \frac{(q^{2(m+j-k+1)},q^{2(\omega-m-j+1)}, q^{-2(N-3-\lambda+k-m-j)-\omega-1}/\overline
z)_k} {(q^{2(\lambda-k+2)}, q^{-2(N-1)+\omega+1}/\overline z)_k}\\
 &\times \frac{(q^{-2(\lambda+j+m+1-k)}, q^{-2(N-k+m-1)+\omega+1}/\overline
z)_{N-j-k}} {(q^{2(\lambda+2-N+j)},
q^{-2(N-j-k)+\omega+1}/\overline z)_{N-j-k}}\\
 &\times {}_{10}\omega_9[q^{2(\lambda-k+1)}; q^{-2k}, q^{-2j}, q^{2(\lambda+j-N+1)},
 q^{2(\lambda-k-\omega+m+j+1)},\\
  & \qquad q^{2(\lambda+j+m-k+2)},   \overline zq^{2(N-m-j)+\omega-1},
 q^{-2(m+j-1)+\omega+1}/\overline
z].
    \end{split}
    \end{equation*}
\end{lemma}

Lemmas \ref{lemma:pit} and \ref{lemma:pit*} can be used to convert
the relations (\ref{eqn:unorth}) to bi-orthogonality relations for
elliptic hypergeometric series. The resulting bi-orthogonality
relations of Theorem \ref{thm:biorth} and \ref{thm:dualbiorth}
have been obtained previously by Frenkel and Turaev
\cite{FrenkelTuraev} and Spiridonov and Zhedanov \cite{SpiriZhe}
(see also Remark \ref{rem:spiri}).

\begin{thm}\label{thm:biorth}
A bi-orthogonality relation for the elliptic hypergeometric series
is given by
    \begin{equation*}
    \begin{split}
\delta_{kl}h_k=& \sum_{j=0}^N w_j \; {}_{10}\omega_9
[q^{2(\Lambda-2l-j+1)};q^{-2j}, q^{-2l},
q^{2(\Lambda-l-N+1)},q^{2(\Lambda-l-\omega+M+1)},\\
     & \qquad q^{2(\Lambda-l+M+2)},  z q^{2(N-M-j-l)+\omega-1},
 q^{-2(M+j+l-2)+\omega-1}/z]\\
     & \times{}_{10}\omega_9[q^{2(\Lambda-2k-j+1)};q^{-2j}, q^{-2k},
q^{2(\Lambda-k-N+1)},q^{2(\Lambda-l-\omega+M+1)},\\
     & \qquad q^{2(\Lambda-k+M+2)}, z q^{2(N-M-j-k)+\omega-3},
 q^{-2(M+j+k-2)+\omega+1}/ z],
    \end{split}
    \end{equation*}
where the quadratic norm $h_k$ and the weight function $w_j$ are
given by
    \begin{equation*}
    \begin{split}
h_k=& \frac{(q^2,q^{-2(\Lambda+M+1)}, q^{-2(\Lambda-\omega+M)},
q^{-2(\Lambda-N)})_k}
    {(q^{2(M+1)}, q^{-2N}, q^{-2(\omega-M)},
    q^{-2\Lambda})_k}\frac{(q^{-2(\Lambda+1)})_{2k}}{(q^{-2\Lambda})_{2k}}\\
 &\times\frac{( zq^{2M-\omega-1},  q^{2(M-N)-\omega+5}/z)_k}
 {( q^{-2(\Lambda+M)+\omega+1}/z,  zq^{2(N-\Lambda-M)+\omega-5})_k}\\
 &\times\frac{(q^{-2(\Lambda-\omega+2M+1)}, q^{-2\Lambda})_N}
    {(q^{-2(\Lambda+M+1)}, q^{-2(\Lambda-\omega+M)})_N}
    \frac{( zq^{\omega-1},z q^{-\omega-3})_N}
    {(  z q^{2M-\omega-1},  z q^{-2M)+\omega-3})_N},
    \end{split}
    \end{equation*}
and $w_j=:w_1(j,k)w_2(j,l)$ with
    \begin{equation*}
    \begin{split}
 w_1(j,k)=& q^{2j-2k} \frac{\theta(q^{2(\Lambda-\omega+2M-N+1+2j)})}{\theta(q^{2(\Lambda-\omega+2M-N+1)})}
   \frac{(q^{2(\Lambda-\omega+2M-N+1)})_j}{(q^{2(\Lambda-\omega+2M+2)})_j}
   \frac{( zq^{2(\Lambda-k+M)-\omega-1})_j}{( q^{-2(N-2-M-k)-\omega+1}/ z)_j}\\
 &\times\frac{(q^{2(M+k+1)}, q^{-2(N-k)}, q^{-2(\Lambda-k+1)}, q^{-2(\omega-M-k)} )_j}
  {(q^2, q^{2(\Lambda-N+M+2)}, q^{2(M+1)}, q^{-2(\Lambda-2k+1)},
   q^{-2N}, q^{-2(\omega-M)}, q^{2(\Lambda-N-\omega+M+1)})_j},\\
 w_2(j,l)=&\frac{(q^{2(M+l+1)},q^{-2(N-l)}, q^{-2(\Lambda-l+1)}, q^{-2(\omega-M-l)})_j}
  {(q^{-2(\Lambda-2l+1)})_j}
   \frac{( q^{-2(N-3-\Lambda+l-M)-\omega-1}/z)_j}
  {( z q^{2(M+l)-\omega-1})_j}.
     \end{split}
     \end{equation*}
\end{thm}

\begin{remark}
These relations are bi-orthogonality relations since there is a
shift in the spectral parameter $z$. Omitting all other parameters
the bi-orthogonality relations are in fact relations of the form
    \begin{equation*}
\delta_{kl}h_k=\sum_{j} w_j P_l(j,q^2z) P_k(j,z).
    \end{equation*}
\end{remark}

\begin{proof}
Applying the dynamical representation $\pi^\omega$ of Proposition
\ref{prop:dynrepr} to (\ref{eqn:unorth1}) gives
    \begin{equation*}
    \begin{split}
\delta_{kl}e_m=&\sum_{j=\max(0,k-m)}^N
\tilde\tau^{N\omega}_{j,l,m+j-k}(\lambda,z)\frac{\Gamma_j(\lambda-\omega-N+2m+2j-2k+2l)}{\Gamma_k(\lambda-N+2l)}\\
 &\times\left[\prod_{i=0}^{N-1}q^{-\omega}\frac{\theta_p( q^{-2i+1+\omega}/\overline z)}
 {\theta_p(q^{-2i+1\omega}/\overline z)}\right]\tau^{N\omega}_{jkm}(\lambda-N+2l, q^{-2(N-2)}/\overline
z)e_{m-k+l}.
    \end{split}
    \end{equation*}
Replacing $\lambda+2l$ by $\Lambda$, $m-k$ by $M$ and $z$ by
$\overline z$ we obtain
    \begin{equation*}
    \begin{split}
\delta_{kl}=& \sum_{j=\max(0,-M)}^N
\frac{\Gamma_j(\Lambda-\omega+2M+2j-N)}{\Gamma_k(\Lambda-N)}\prod_{i=0}^N
   q^{-\omega}\frac{\theta_p(q^{-2i+1+\omega}/ z)}
   {\theta_p(q^{-2i+1\omega}/ z)}\\
 &\times \tilde\tau^{N\omega}_{j,l,M+j}(\Lambda-2l,\overline z)
    \tau^{N\omega}_{j,k,M+k}(\Lambda-N,q^{-2(N-2)}/z).
    \end{split}
    \end{equation*}
Using Lemmas \ref{lemma:pit} and \ref{lemma:pit*} and elementary
relations for the elliptic factorials proves the theorem.
\end{proof}

\begin{thm}\label{thm:dualbiorth}
The dual bi-orthogonality relation for the elliptic hypergeometric
series is given by
    \begin{equation*}
    \begin{split}
\delta_{kl}=&\sum_{j} \frac{w_1(l,j)w_2(k,j)}{(h_j)}\;
{}_{10}\omega_9[q^{2(\Lambda-2j-l+1)}; q^{-2l}, q^{-2j},
q^{2(\Lambda-N-j+1)}, q^{2(\Lambda-j-\omega+M+1)},\\
 & \qquad q^{2(\Lambda+M-j+2)}, q^{-2(M+j+l-2)+\omega+1}/ z, z q^{-2(M+j-N+l+1)+\omega-1}]\\
 &\times{}_{10}\omega_9[q^{2(\Lambda-2j-k+1)};
q^{-2k}, q^{-2j}, q^{2(\Lambda-N-j+1)},
q^{2(\Lambda-j-\omega+M+1)},\\
 & \qquad q^{2(\Lambda+M-j+2)}, q^{-2(M+j+k-1)+\omega+1}/z, z q^{-2(M+j-N+k)+\omega-1}],
    \end{split}
    \end{equation*}
where $w_1$, $w_2$ and $h_j$ are as in Theorem \ref{thm:biorth}.
\end{thm}

\begin{proof}
These dual bi-orthogonality relations can be computed from
(\ref{eqn:unorth2}) by applying the dynamical representation.
Since the biorthogonal system in Theorem \ref{thm:biorth} is known
to be self-dual \cite{SpiriZhe}, we can also obtain the dual
relations from Theorem \ref{thm:biorth}.
\end{proof}

\begin{remark}\label{rem:spiri}
In \cite{FrenkelTuraev} an elliptic analogue of Bailey's
transformation formula is proved. Let $bcdefg=a^3q^{2(n+2)}$ and
$\lambda=a^2q^2/bcd$. Then
    \begin{equation}\label{eqn:bailey}
    \begin{split}
{}_{10}\omega_9&[a;b,c,d,e,f,g,q^{-2n}]=\frac{(aq^2,aq^2/ef,\lambda
q^2/e, \lambda q^2/f)_n}{(aq^2/e,aq^2/f,\lambda q^2/ef, \lambda
q^2 )_n}{}_{10}\omega_9[\lambda;\lambda b/a, \lambda c/a, \lambda
d/a,e,f,g,q^{-2n}]
    \end{split}
    \end{equation}
We can relate the bi-orthogonality relations of Theorem
\ref{thm:biorth} and \ref{thm:dualbiorth} to the ones given in
\cite{SpiriZhe}. To obtain this relation explicitly we have to
apply the elliptic analogue of Bailey's transformation formula
(\ref{eqn:bailey}) twice to both $_{10}\omega_9$-functions in our
bi-orthogonality relations in different ways. Finally, let us
emphasize that we do not need Bailey's transformation formula to
obtain the bi-orthogonality relations of Theorem \ref{thm:biorth}
and \ref{thm:dualbiorth} in the symmetric form given.
\end{remark}

\begin{remark}
Using the dynamical representation of Proposition
\ref{prop:dynrepr} we can obtain transformation formula
\eqref{eqn:bailey} from the unitarity property of the
corepresentations stated in Proposition \ref{prop:unitarity}.
\end{remark}



\end{document}